\documentclass [12 pt]{amsart}
\usepackage{amsmath,amssymb,amsthm,hyperref}
\usepackage{a4wide}
\usepackage{color}
\usepackage{enumerate}
\usepackage{cleveref}
\usepackage{tikz}
\usepackage{tikz}
\tikzset{font={\fontsize{1pt}{12}\selectfont}}
\usetikzlibrary{decorations.pathreplacing}
\usetikzlibrary{shapes,snakes}
 \newcounter{case}
 
 \renewcommand{\thecase}{\arabic{case}}
\newtheorem{thm}{Theorem}[section]
\newtheorem*{thm*}{Theorem}
\newtheorem*{conj*}{Conjecture}

\newtheorem{fact}[thm]{Fact}

\theoremstyle{remark}

\theoremstyle{definition}
\newtheorem{eg}[thm]{Example}
\newtheorem{defn}[thm]{Definition}

\newcounter{claim}[thm]

\newcounter{subcase}

\renewcommand{\thesubcase}{\Alph{subcase}}

\title {On the singular pebbling number of a graph}

\author {Harmony Morris}

\email {harmony.morris05@gmail.com}
\begin {document}

\begin {abstract}
In this paper, we define a new parameter of a graph as a spin-off of the pebbling number (which is the smallest $t$ such that every supply of $t$ pebbles can satisfy every demand of one pebble). This new parameter is the singular pebbling number, the smallest $t$ such that a player can be given any configuration of at least $t$ pebbles and any target vertex and can successfully move pebbles so that exactly one pebble ends on the target vertex. We also prove that the singular pebbling number of any graph on 3 or more vertices is equal to its pebbling number, that the singular pebbling number of the disconnected graph on two vertices is equal to its pebbling number, and we find the singular pebbling numbers of the two remaining graphs, $K_1$ and $K_2$, which are not equal to their pebbling numbers.
\end {abstract}

\maketitle

\section {Introduction}

In this paper, all graphs are simple (loops and multiple edges are not allowed). A complete graph is denoted by $K_n$, with $n$ being the number of vertices. The minimum valency of a graph $G$ is denoted $\delta(G)$.

Graph pebbling is a game in which a player moves pebbles from vertices to other vertices along edges in an attempt to end up with at least one pebble on the `target vertex'. The player has an opponent who decides the initial placement of the pebbles and chooses a vertex to be the target vertex. The player can move pebbles by taking two pebbles off of any vertex with two or more pebbles and moving one of them to an adjacent vertex, removing the other pebble from the game.

The problem of graph pebbling was initially studied in \cite{chung}. Two survey papers have been written on the various research done on this problem, \cite{hurlbert2004} and \cite{hurlbert2013}.

Now, we define the pebbling number of a graph $G$.

\begin {defn}[\cite{hurlbert2011}]
A graph's \emph{pebbling number}, denoted $\pi(G)$, is the smallest $t$ such that every supply of $t$ pebbles can satisfy every demand of one pebble. 
\end {defn}

In other words, the pebbling number is the smallest $t$ such that a player can be given any configuration of $t$ pebbles and any target vertex and can successfully move pebbles so that at least one pebble ends on the target vertex.

This paper introduces a new concept, the singular pebbling number of a graph $G$.

\begin {defn}
A graph's \emph{singular pebbling number}, denoted $\pi_s(G)$, is the smallest $t$ such that a player can be given any configuration of at least $t$ pebbles and any target vertex and can successfully move pebbles so that exactly one pebble ends on the target vertex. 
\end {defn}

This paper also uses a different interpretation of pebbling number: the pebbling number of a graph is one more than the largest $t$ such that a player can configure $t$ pebbles and choose a target vertex so that another player cannot move pebbles to get at least one pebble on the target vertex. Thus, the singular pebbling number of a graph is one more than the largest $t$ such that a player can configure $t$ pebbles and choose a target vertex so that another player cannot move pebbles to get exactly one pebble on the target vertex. For the remainder of this paper, the `player' will therefore be the person attempting to block the `opponent' from pebbling the graph.

In this paper, we thus assume the goal is for the player to configure $t$ pebbles and choose a target vertex so that their opponent cannot move pebbles to get exactly one pebble on the target vertex.

\section {Preliminary Results}

\begin{fact}\label{fact:targetcan'tstartw4}

The target vertex cannot begin with 3 or more pebbles unless it is an isolated vertex.

\end{fact}

\begin{proof}

If there is an even number ($\ge 4$) of pebbles on the target vertex, all the pebbles can be moved to an adjacent vertex (which means the number of pebbles is divided by two, but since there were at least four pebbles to begin with at least two pebbles remain). There are therefore two pebbles on a vertex adjacent to the target vertex that can be moved to the target vertex, leaving precisely one pebble on the target vertex and allowing the opponent to win the game.

If there is an odd number ($\ge 3$) of pebbles on the target vertex, all but one of the pebbles (an even number) can be moved to an adjacent vertex, leaving exactly one pebble on the target vertex.

Thus, as long as there is a vertex adjacent to the target vertex, the target vertex cannot begin with 3 or more pebbles.

\end{proof}

\begin{fact}\label{fact:implies}
$\pi_s(G) \neq \pi(G)$ implies that the player's optimal configuration of pebbles on $G$ must involve exactly 2 pebbles being placed on the target vertex, unless the target vertex is an isolated vertex.
\end{fact}

\begin{proof}
By Fact \ref{fact:targetcan'tstartw4}, three or more pebbles cannot begin on the target vertex, and if one pebble begins on the target vertex the opponent automatically wins. If there are no pebbles beginning on the target vertex, the opponent will never be able to move more than one pebble onto the target vertex without being able to move exactly one pebble onto the target vertex. This is true because if there are more than two pebbles on a vertex adjacent to the target vertex, the opponent can choose to move just two pebbles toward the target vertex. Thus, if an optimal configuration can be achieved without any pebbles being placed on the target vertex, $\pi_s(G) = \pi(G)$.
\end{proof}

\begin{fact}\label{fact:pebsontarget}
2 pebbles being on the target vertex means no pebbles can ever be moved onto any vertex adjacent to the target vertex (except off of the target vertex) or the player loses.
\end{fact}

\begin{proof}
If there is ever a pebble on a vertex adjacent to the target vertex, the two pebbles on the target vertex can be moved off the target vertex toward that vertex. There are then two pebbles on a vertex adjacent to the target vertex, which can subsequently be moved to leave exactly one pebble on the target vertex.
\end{proof}

\section {Main Results}

\begin {thm} \label{thm:piseqpi}
For all graphs on three or more vertices and for the disconnected graph on two vertices, $\pi_s(G) = \pi(G)$.
\end {thm}

\begin {proof}
For any disconnected graph on two or more vertices, $\pi_s(G) = \pi(G)$ because for both, an infinite number of pebbles can simply be placed on a vertex that is not connected to the target vertex, making both the pebbling number and singular pebbling number infinite.

We can thus assume that for $\pi_s(G) \neq \pi(G)$, $\delta(G) > 0$.

Additionally, if $\delta(G) \ge 2$, $\pi_s(G) = \pi(G)$. By Fact \ref{fact:implies}, the only way $\pi_s(G)$ can be different from $\pi(G)$ on a connected graph with two or more vertices is if the player's optimal configuration of pebbles must involve 2 pebbles being placed on the target vertex. By Fact \ref{fact:pebsontarget}, 2 pebbles being on the target vertex means no pebbles can ever be moved onto the vertices adjacent to the target vertex. Given any pebble configuration that prevents the opponent from getting a single pebble to the target vertex after starting with two pebbles on the target vertex, one can simply add a pebble to each vertex adjacent to the target vertex (which must be at least two vertices) and remove the pebbles from the target vertex, allowing pebbles on the vertices adjacent to the target vertex but still not allowing a pebble to reach the target vertex. Thus, placing two pebbles on the target vertex is not necessary to achieve an optimal pebble placement. If every vertex of $G$ has valency two or more, since the target vertex is a vertex, it must have valency two or more and therefore $\pi_s(G) = \pi(G)$.

This means we can also assume that for $\pi_s(G) \neq \pi(G)$, $\delta(G) < 2$. Thus, since $\delta(G) < 2$ and $\delta(G) > 0$, $\delta(G) = 1$ and there therefore must be at least one vertex with valency exactly 1 in $G$ for $\pi_s(G) \neq \pi(G)$, if $G$ is a graph on three or more vertices. The target vertex must be a vertex of valency 1 given the proof in the paragraph above.

By Fact \ref{fact:implies}, for the singular pebbling number of a connected graph on two or more vertices to be different from the pebbling number, the optimal pebble configuration must have two pebbles on the target vertex. Thus, to prove that $\pi_s(G) = \pi(G)$ it suffices to show that it is equivalent or better not to have any pebbles on the target vertex. 

By Fact \ref{fact:pebsontarget}, 2 pebbles being on the target vertex means no pebbles can ever be moved onto the vertex adjacent to the target vertex. In a graph on three or more vertices in which the target vertex has valency 1, there is at least one vertex at distance two or more from the target vertex. Given any pebble configuration that prevents the opponent from getting a single pebble to the target vertex after starting with two pebbles on the target vertex, one can simply add two pebbles to a vertex at distance 2 from the target vertex and remove the pebbles from the target vertex, allowing a pebble to get to the vertex adjacent to the target vertex but still not allowing a pebble to get to the target vertex. Thus, we have shown that it is always equivalent or better for the player not to place 2 pebbles on the target vertex for connected graphs on at least three vertices with $\delta(G) = 1$ and we have therefore concluded the proof of this theorem.

\end {proof}

\begin {eg} \label{eg:k1}
$\pi(K_1) = 1$, whereas $\pi_s(K_1)$ is infinite.
\end {eg}

\begin {proof}
If there are 0 pebbles on any graph, it is impossible for the opponent to move a pebble to the target vertex. If there is one pebble on a graph with only one vertex, it is impossible not to have at least one pebble on the target vertex. The pebbling number of $K_1$ is thus one.

The player can put as many pebbles as they want (at least two) onto the vertex of $K_1$ in finding the singular pebbling number, since there can never be exactly one pebble on the vertex if the vertex begins with more than one pebble. Thus, the singular pebbling number of $K_1$ is infinite.
\end {proof}

\begin {eg} \label{eg:k2}
$\pi(K_2) = 2$, whereas $\pi_s(K_2) = 3$.
\end {eg}

\begin {proof}
The player can only place pebbles on the vertex adjacent to the target vertex in finding the pebbling number. If the player places more than one pebble on the adjacent vertex, the opponent will win by moving two pebbles toward the target vertex. The pebbling number of $K_2$ is thus two.

The player can either place one pebble on the vertex adjacent to the target vertex or can place two pebbles on the target vertex in finding the singular pebbling number. The singular pebbling number of $K_2$ is thus three.
\end {proof}

\textsc{Acknowledgements.} The author thanks Joy Morris for her help with citations, as well as her explanations of previous research on pebbling number and her helpful suggestions throughout the process. The author also thanks Dave Morris for his helpful suggestions throughout the process.

\end {document}